\documentclass[12pt]{amsart}
\usepackage{a4wide}
\usepackage{amssymb}
\newenvironment{pf}{\medskip\noindent{\it Proof.}\enspace}%
  {\hfill$\square$\newline\smallskip}
\newtheorem{lem}{Lemma}
\newtheorem{theo}[lem]{Theorem}

\newtheorem{prop}[lem]{Proposition}
\newtheorem{rmk}[lem]{Remark}

\def\q{\hskip0.17cm}
\def\,{\hskip0.12cm}

\def\frak{\mathfrak}
\def\Bbb{\mathbb}
\def\cal{\mathcal}

\def\m{\begin{pmatrix}}
\def\em{\end{pmatrix}}
\def\sm{\left(\smallmatrix}
\def\esm{\endsmallmatrix\right)}

\begin{document}
\voffset=24pt
\title{\bf More on super-replication formulae}

\author{Chang Heon Kim \and Ja Kyung Koo}

\address{Chang Heon Kim, Department of mathematics, Seoul Women's university,
126 Kongnung 2-dong, Nowon-gu, Seoul, 139-774 Korea}
\email{chkim@swu.ac.kr}

\address{Korea Advanced Institute of Science and
Technology, Department of Mathematics, Taejon, 305-701 Korea}
\email{jkkoo@math.kaist.ac.kr}

 \begin{abstract}
 We extend Norton-Borcherds-Koike's replication formulae
 to super-replicable ones by working with the congruence groups $\Gamma_1(N)$
 and find the product identities which characterize
 super-replicable functions. These will provide a clue for
 constructing certain new infinite dimensional Lie superalgebras whose
 denominator identities coincide with the above product identities. Therefore it
 could be one way to find a connection between modular functions
 and infinite dimensional Lie algebras.
 \end{abstract}

 \maketitle
 \renewcommand{\thefootnote}%
             {}
 \footnotetext{Supported by KOSEF Research Grant 98-0701-01-01-3
 \par AMS
  Classification : 11F03,
       11F22}
 \baselineskip=24pt

\section{Introduction}
 \par
 Let ${\Bbb M}$ be the monster simple group whose order is
 approximately $8\times 10^{53}$.
 In 1979, Mckay and Thompson (\cite{M-T}) found some relations between the monster group
 and elliptic modular function $J$. Let $f_i$ be the degrees of
 irreducible characters of ${\Bbb M}$ (i.e. $f_1=1$, $f_2=196883$,
 $f_3=21493760$, etc) and $J=j-744=q^{-1}+0+196884 q +
 \cdots = q^{-1}+0+\sum_{n\ge 1} c_n q^n$ be the elliptic modular
 function (called the modular invariant). Then it was noticed that $c_i$ can be
 expressed as a linear
 combination of $f_j$'s. For example,
 $c_1=f_1+f_2$, $c_2=f_1+f_2+f_3$,
 $c_3=2f_1+2f_2+f_3+f_4$, etc.
 Based on their observation, Mckay and Thompson conjectured the
 existence of a natural infinite dimensional representation of the
 monster, $V=\bigoplus_{n\in {\Bbb Z}} V_n$
 with ${\rm dim}(V_n)=c_n$.
 Thompson (\cite{Thompson}) also proposed considering, for any
 element
 $g\in {\Bbb M}$, the modular properties of the series
 $$ T_g(q)=q^{-1} + 0 + \sum_{n\ge 1} {\rm Tr} (g|V_n) q^n $$
 (the {\it Thompson series})
 where $q=e^{2\pi i z}$ with
 Im$(z)>0$. And Conway and
 Norton (\cite{C-N}) conjectured that
 $T_g(q)$ would be a Hauptmodul for a suitable genus zero
 discrete subgroup of $SL_2({\Bbb R})$.
 \par Now we shall explain the notion of a Hauptmodul.
 Let ${\frak H}$ be the complex upper half plane and let $\Gamma$
 be a discrete subgroup of $SL_2({\Bbb R})$. Since
 the group $\Gamma$ acts on ${\frak H}$ by linear fractional
 transformations, we get the modular curve
 $X(\Gamma)=\Gamma\backslash{\frak H}^*$, as a projective closure
 of the smooth affine curve $\Gamma\backslash{\frak H}$.
 Here, ${\frak H}^*$ denotes the union of ${\frak H}$ and $\Bbb
 P^1({\Bbb Q})$. If the genus of $X(\Gamma)$ is zero, then
 the function field $K(X(\Gamma))$ in this case is a
 rational function field ${\Bbb C}(j_{\Gamma})$ for some modular
 function $j_{\Gamma}$. If the element $\sm 1&1\\0&1 \esm$ is in
 $\Gamma$, it takes
 $z$ to $z+1$, and in particular a modular function $f$ in
 $K(X(\Gamma))$ is periodic. Thus it can be written as a Laurent series in
 $q=e^{2\pi i z}$ ($z\in {\frak H}$), which is called a $q$-{\it series} \, (or
 $q$-{\it expansion}) \, of $f$. We call $f$ {\it normalized} \, if
 its $q$-series starts with $q^{-1}+0+a_1 q+a_2 q^2+\cdots$. By a
 {\it Hauptmodul} $t$ we mean the normalized generator of
 a genus zero function field $K(X(\Gamma))$.

 \par In 1980s Frenkel, Lepowsky and Meurman (\cite{Frenkel1}, \cite{Frenkel2})
 constructed an infinite
 dimensional graded representation
 $V=\bigoplus_{n\in {\Bbb Z}} V_n$ of ${\Bbb M}$ with dim$(V_n)=c_n$ using vertex operator
 algebra, which comes from
 the space of physical states of a bosonic string
 moving in a ${\Bbb Z}_2$-orbifold of a 26-dimensional torus.
 In 1992 Borcherds (\cite{Borcherds}) proved the Conway-Norton's conjecture by using a
 recursion formula derived from the twisted denominator formula of
 a generalized Kac-Moody Lie algebra with some group action.
 Here, we shall investigate some properties of the twisted denominator
 formula.
 \par Put $t=T_g$ and $t^{(i)}=T_{g^i}$. Let
 $H_n={\rm Tr} (g|V_n)$ and $H_n^{(i)}={\rm Tr} (g^i|V_n)$. Also
 let $X_n(t)$ be the Faber polynomial of $t$ such that $X_n(t)\equiv \frac 1n
 q^{-n} \mod q{\Bbb C}[[q]]$ (\cite{Duren}, Chapter 4)
 and write $X_n(t)=\frac 1n q^{-n}+\sum_{m\ge 1} H_{m,n} q^m$.
 Then the Borcherds' twisted denominator formula is
 $$ p^{-1}\prod_{n>0\atop m\in {\Bbb Z}} \exp(-\sum_{s\ge 1} \frac 1s H_{mn}^{(s)}
    q^{sm} p^{sn}) = t(p)-t(q) $$ where $p=e^{2\pi i y}$ and $q=e^{2\pi i
    z}$ with $y,z\in {\frak H}$.
 Moreover it is not difficult to show that
 the above is equivalent to $$X_n(t)=n^{-1}\sum_{ad=n \atop 0\le b
 < d} t^{(a)}\left(\frac{az+b}{d}\right).$$
 By comparing $q^m$-terms on both sides we obtain
 $$ H_{m,n}=\sum_{s|(m,n)} H_{mn/d^2}^{(s)} \, ,$$ which is equivalent
 to $$H_{a,b}=H_{c,d} \text{ whenever } (a,b)=(c,d) \text{ and }
 ab=cd \text{ (see \cite{Koike})}.$$
 When $t=\frac 1q+\sum_{m\ge 1} H_m q^m$ is a formal power series with
 integer coefficients, we call
 $t$ {\it replicable} if it satisfies
 $H_{a,b}=H_{c,d} \text{ whenever } (a,b)=(c,d) \text{ and }
 ab=cd$. One then has the following equivalent conditions which
 characterize replicable functions
 (see \cite{Kang}, also \cite{ACMS}, \cite{Koike} and \cite{Norton} for (a)$\iff$(b)).
 \par (a) $t$ is replicable.
 \par (b) for each $s\ge 1$ there exists a normalized $q$-series $t^{(s)}=\frac
 1q+\sum_{m\ge 1} H_m^{(s)} q^m$ such that
 $$ t^{(1)}=t, \text{ \q and \q } X_n(t)=n^{-1}\sum_{ad=n \atop 0\le b < d}
 t^{(a)}\left(\frac{az+b}{d}\right).$$
 \par (c) for every $s\ge 1$ there exists a normalized $q$-series $t^{(s)}=\frac
 1q+\sum_{m\ge 1} H_m^{(s)} q^m$ such that
 $$ p^{-1}\prod_{n>0\atop m\in {\Bbb Z}} \exp(-\sum_{s\ge 1} \frac 1s H_{mn}^{(s)}
    q^{sm} p^{sn}) = t(p)-t(q) $$ where $p=e^{2\pi i y}$ and $q=e^{2\pi i
    z}$ with $y,z\in {\frak H}$.
 \par\noindent
  We note that when $t$ is the Thompson series, the product identity (c) coincides with
 the Borcherds' twisted denominator formula.

 \par Let $t_0$ be a replicable function. According to Norton's conjecture (\cite{Norton})
 $t_0$ is a Hauptmodul for some modular curve $X$.
 Let $\tilde{X}$ be a double covering of $X$ with genus zero and
 $t$ be the Hauptmodul for $\tilde{X}$.
 It is natural to ask what sort of relations $t_0$ and $t$ satisfy.
 As before we write $t_0=\frac 1q+\sum_{m\ge 1} h_m q^m$ (resp.
 $t=\frac
 1q+\sum_{m\ge 1} H_m q^m$) and
 $X_n(t_0)=\frac 1n q^{-n}+\sum_{m\ge 1} h_{m,n} q^m$ (resp.
 $X_n(t)=\frac 1n q^{-n}+\sum_{m\ge 1} H_{m,n} q^m$).
 We propose to call $t$ {\it super-replicable} with respect to a replicable
 function $t_0$ and a character $\psi$ if the following formula
 holds
 $$ \psi (a,b) \times (2H_{a,b}-h_{a,b}) =
 \psi (c,d) \times (2H_{c,d}-h_{c,d})$$
 whenever $(a,b)=(c,d)$ and $ab=cd$.
 Observe that if $H_{a,b}=h_{a,b}$ and $\psi$ is the trivial character, then the
 above identity reduces to the replicable one.
 \par
 Let $\Gamma_1(N)=\{\sm a & b \\ c& d \esm \in SL_2({\Bbb Z}) \mid
         c\equiv 0 \mod N, a\equiv d\equiv 1 \mod N \} $
 and
  $X_1(N)= \Gamma_1(N)\backslash {\frak H}^*$, and let
 $\Gamma_0(N)=\{ \sm
   a & b \\ c& d \esm \in SL_2({\Bbb Z}) \mid
         c\equiv 0 \mod N \} $
 and
  $X_0(N)= \Gamma_0(N)\backslash {\frak H}^*$.
 According to \cite{Kim-Koo1}, \cite{Miyake} the genus of
 $X_1(N)$ is
 zero if and only if $1\le N \le 10$ and $N=12$.
 When $N\in \{5,8,10,12 \}$, $X_1(N)$ becomes a double covering of
 $X_0(N)$.
 Let $t$ (resp. $t_0$) be the Hauptmodul of $X_1(N)$ (resp.
 $X_0(N)$) for such $N$.
 We proved in \cite{Super}, Corollary 11 that
 for positive integers $a,b,c,d$ with
 $ab=cd$, $(a,b)=(c,d)$ and $(b,N)=(d,N)=1$,
 \begin{equation}
 \psi_N (b) \times (2H_{a,b}-h_{a,b}) =
 \psi_N (d) \times (2H_{c,d}-h_{c,d})
 \tag{$*$}
 \end{equation}
 where $\psi_N: \, ({\Bbb Z}/N{\Bbb Z})^\times \to \{ \pm 1\}$ is a
 character defined by
 $$ \psi_N(b)= \begin{cases} 1, & \text{ if } b\equiv \pm 1 \mod N \\
 -1, & \text{otherwise.} \end{cases} $$
 In this article we refine the above to discard the condition $(b,N)=(d,N)=1$ and show
 that each Hauptmodul $t$ is in fact a super-replicable function with respect to a
 replicable function $t_0$ and some charcter $\psi$
 (Theorem \ref{BB}, \ref{GG} and \ref{KK}).
 And, in Theorem \ref{product2} we give a criterion which characterizes
 super-replicable functions.

 \par Through the article we adopt the following notations:
 \par\noindent $\bullet$ \, $f|_{\sm a&b \\ c&d
 \esm}=f\left(\frac{az+b}{cz+d}\right)$
 \par\noindent $\bullet$ \,
     $f(z)=g(z)+O(1) \q \text{ means that $f(z)-g(z)$ is bounded as $z$
                    goes to $i\infty$. }$

 \section{Super-replication formulae}
 \par
 Let $\Delta^n$ be the set of $2\times 2$ integral matrices
 $\sm a&b \\ c&d \esm$ where $a\in 1+N{\Bbb Z}, \, c\in N{\Bbb Z}$
 and $ad-bc=n$. Then $\Delta^n$ has the following right coset
 decomposition: (See \cite{Koblitz}, \cite{Miyake},
 \cite{Shimura1})
 $$
 \Delta^n=\bigcup_{\smallmatrix a|n \\ (a,N)=1 \endsmallmatrix}
     \bigcup_{i=0}^{\frac{n}{a}-1} \Gamma_1(N)\sigma_a
          \m a&i \\ 0&\frac{n}{a} \em
 $$
 where $\sigma_a\in SL_2({\Bbb Z})$ such that $\sigma_a\equiv \sm
 a^{-1}&0 \\ 0&a \esm \mod N$.
 For a modular function $f$ we define the Hecke operators
 $U_n$ and $T_n$ by
 \begin{align*}
 f|_{U_n} &=n^{-1} \sum_{i=0}^{n-1} f|_{\sm 1&i\\0&n \esm} \\
 \intertext{and}
 f|_{T_n} &= n^{-1} \sum_{\smallmatrix a\mid n \\ (a,N)=1
 \endsmallmatrix}
   \sum_{i=0}^{\frac na-1} f|_{\sigma_a \sm a& i \\ 0& \frac na \esm}.
 \end{align*}
 Notice that if $n$ divides a power of $N$, then
 $f|_{U_n}=f|_{T_n}$.
 \begin{lem}
 If $f$ is invariant under $\Gamma_1(N)$ and $\gamma_0\in
 \Gamma_0(N)$, then
 $(f|_{T_n})|_{\gamma_0}=(f|_{\gamma_0})|_{T_n}$ for any positive
 integer $n$. In particular, $f|_{T_n}$ is again invariant under $\Gamma_1(N)$.
 \label{Hecke}
 \end{lem}
 \begin{pf}
 \cite{Super}, Lemma 7.
 \end{pf}
  Let $t$ (resp. $t_0$) be the Hauptmodul of $\Gamma_1(N)$ (resp.
 $\Gamma_0(N)$).
 \begin{lem}
 For $\sm a&b\\c&d \esm \in \Gamma_0(N)$, we have
 $$ X_r(t)|_{\sm a&b\\ c&d \esm}
  = \frac 12 \{\psi_{N}
 (d)(2X_r(t)-X_r(t_0))+X_r(t_0)\}+c_0 $$
 with $c_0=
   \begin{cases} 0, & \text{ if } d\equiv \pm 1 \mod 12 \\
                 X_r(t)(a/c), & \text{ otherwise. }
   \end{cases}$
 \label{Action}
 \end{lem}
 \begin{pf}
 Observe that
 $$ X_r(t)|_{\sm a&b\\c&d \esm}+X_r(t)=
 \begin{cases} 2 X_r(t), &  \text{ if } d\equiv \pm 1 \mod N \\
               X_r(t_0)+ X_r(t)(a/c), &  \text{ otherwise.}
 \end{cases} $$
 In the above, when $d$ is not congruent to $\pm 1 \mod N$, $
 X_r(t)|_{\sm a&b\\c&d \esm}+X_r(t)$ is invariant under $\Gamma_0(N)$ and has poles only
 at $\Gamma_0(N)\infty$ with $r^{-1}q^{-r}$ as its pole part. This
 guarantees the above equality. Now the assertion is immediate.
 \end{pf}
 \begin{lem}
 Let $N\in \{5,8,10,12 \}$ and $p$ be a prime dividing $N$.
 Except the case $N=10$ and $p=2$,
 we have \, \,
 $2X_{pn}(t)|_{U_p} - X_{pn}(t_0)|_{U_p}=\frac 1p
 (2X_n(t)-X_n(t_0))$.
 \label{AA}
 \end{lem}
 \begin{pf}
 First we observe that all functions in the assertion are
 invariant under $\Gamma_1(N)$. By definitions,
 $X_n(t)$ (resp. $X_n(t_0)$) can have poles only at $\Gamma_1(N)\infty$
 (resp. $\Gamma_0(N)\infty$) and
 $X_{pn}(t)|_{U_p}$ (resp. $X_{pn}(t_0)|_{U_p}$)
 can have poles only at $\sm 1&i \\ 0&p \esm^{-1}\Gamma_1(N)\infty$
 (resp. $\sm 1&i \\ 0&p \esm^{-1}\Gamma_0(N)\infty$) for
 $i=0,\cdots,p-1$.
 On the other hand, we have
 \begin{align*}
 \sm 1&i \\ 0&p \esm^{-1} \Gamma_0(N)\infty
  &=p^{-1}\sm p&-i \\ 0&1 \esm\Gamma_0(N)\infty
    \sim \sm 1&-i \\ 0&1 \esm
       \sm p&0 \\ 0&1 \esm\Gamma_0(N)\infty \\
  & \sim \sm p&0 \\ 0&1 \esm\Gamma_0(N)\infty
 \end{align*}
 where the symbol $\sim$ stands for a $\Gamma_1(N)$-equivalence of
 cusps.
 Let $\sm a&b\\c&d \esm$ be an element in $\Gamma_0(N)$.
 Then $\sm p&0 \\ 0&1 \esm \sm a&b\\c&d \esm \infty
      =\frac{a}{c/p}$ in lowest terms. Write $c=N\cdot(pk+l)$
 with $k\in\Bbb Z$ and $l\in\{0,\cdots,p-1\}$. We then readily see
 that $\frac cp\equiv\frac Np\cdot l \mod N$, which implies $\sm a
 \\ c/p \esm \equiv \sm a \\ (N/p)l \esm\mod N$.
 \par\noindent
 We claim that for $l\ge 1$, $\frac{a}{(N/p)\cdot l} \sim
 \frac{1}{(N/p)\cdot j}$ for some $j$ in $X_1(N)$.
 In fact, $\frac{a}{(N/p)\cdot l} \sim
 \frac{1}{(N/p)\cdot j}$ for some $j$ in $X_1(N)$
 $\iff \pm \sm a \\ (N/p)l \esm \equiv \sm 1+(N/p)jn \\ (N/p)j \esm\mod
 N$ for some $n$
 $\iff \pm a \equiv 1+(N/p)jn \mod N$ where $j=l$ or $p-l$
 $\iff ((N/p)j,N) | \pm a -1$
 $\iff (N/p) | a-1 \text{ or } a+1$ for $a$ prime to $N$, which is
 clear from the following table:
 \vskip0.3cm
 \begin{center}
 \begin{tabular}{|c|c|c|c|}
 \hline
 $N$ & $p$ & $N/p$ & $a$ \\
 \hline
 $5$ & $5$ & $1$ & $1,2,3,4$ \\
 \hline
 $8$ & $2$ & $4$ & $1,3,5,7$ \\
 \hline
 $10$ & $5$ & $2$ & $1,3,7,9$ \\
 \hline
 $12$ & $2$ & $6$ & $1,5,7,11$ \\
 \hline
 $12$ & $3$ & $4$ & $1,5,7,11$ \\
 \hline
 \end{tabular}
 \end{center}
 \vskip0.3cm
 Now we consider the pole part of $X_{pn}(t)|_{U_p}$ at the cusp
 $\frac{1}{(N/p)j}$ for $j=1,\cdots,p-1$.
 \begin{align}
 p\cdot(X_{pn}(t)|_{U_p})|_{\sm 1&0 \\ (N/p)j &1 \esm}
 & = \sum_{i=0}^{p-1} X_{pn}(t)|_{\sm 1&i \\0&p \esm \sm 1&0 \\ (N/p)j&1 \esm}
   = X_{pn}(t)|_{\sm 1&0 \\ Nj &p \esm}+\sum_{i=1}^{p-1}
     X_{pn}(t)|_{\sm 1+(N/p)ji&i \\ Nj&p \esm} \notag \\
 & =(X_{pn}(t)|_{\sm 1&0 \\ Nj &1 \esm})|_{\sm 1&0\\0&p \esm}
     +\sum_{i=1}^{p-1} X_{pn}(t)|_{\sm 1+(N/p)ji&i \\ Nj&p \esm} \notag \\
 & = X_{pn}(t)\left(\frac zp \right)
     +\sum_{i=1}^{p-1} X_{pn}(t)|_{\sm 1+(N/p)ji&i \\ Nj&p \esm} \label{pole1}.
 \end{align}
 Similarly
 \begin{align}
 p\cdot(X_{pn}(t_0)|_{U_p})|_{\sm 1&0 \\ (N/p)j &1 \esm}
 & = X_{pn}(t_0)\left(\frac zp \right)
     +\sum_{i=1}^{p-1} X_{pn}(t_0)|_{\sm 1+(N/p)ji&i \\ Nj&p \esm} \label{pole2}.
 \end{align}
 Here we consider the matrix $M_i=\sm 1+(N/p)ji&i \\ Nj&p \esm$
 where $i$ runs through $1,\cdots,p-1$. Write it as $\gamma_i W_i$
 where $\gamma_i\in SL_2({\Bbb Z})$ and $W_i$ is an upper
 triangular matrix. Since $\det \gamma_i W_i=p$, $W_i$ must be of
 the form $\sm p&*\\0&1\esm$ or $\sm 1&*\\0&p\esm$.
 $W_i=\sm p&*\\0&1\esm$ $\iff$ $p$ divides both $1+(N/p)ji$ and $Nj$
 $\iff$ $p$ divides $1+(N/p)ji$.
 $W_i=\sm 1&*\\0&p\esm$ $\iff$ $p$ does not divide $1+(N/p)ji$
 $\iff$ $\gamma_i \in \Gamma_0(N)$.
 \par\noindent
 Case I. $N=p=5$
 \par
 $\{ ji \mod 5 | i=1,2,3,4 \}
 =\{ 1,2,3,4 \mod 5 \} $
 \par
 $\{ 1+(N/p)ji \mod 5 | i=1,2,3,4 \}
 =\{ 2,3,4,0 \mod 5 \} $
 \par
 Three of $\gamma_i$'s belong to $\Gamma_0(5)$ and
 one of them belongs to $\pm\Gamma_1(5)$.
 Thus the formula
 $(\ref{pole1})
  = \frac{q^{-n}}{pn} \times 2 + O(1)$ and the formula
 $(\ref{pole2})
  = \frac{q^{-n}}{pn} \times 4 + O(1)$.
 \par\noindent
 Case II. $N=8$ and $p=2$
 \par
 $i=j=1$ and
 $\{ 1+(N/p)ji \mod 2 \}
 =\{ 5 \mod 2 \} $
 \par
 $\gamma_i$ belongs to $\Gamma_0(8)\backslash \pm\Gamma_1(8)$.
 Therefore
 $(\ref{pole1})
  = \frac{q^{-n}}{pn} + O(1)$ and
 $(\ref{pole2})
  = \frac{q^{-n}}{pn} \times 2 + O(1)$.
 \par\noindent
 Case III. $N=10$ and $p=5$
 \par
 $\{ ji \mod 5 | i=1,2,3,4 \}
 =\{ 1,2,3,4 \mod 5 \} $
 \par
 $\{ 1+(N/p)ji \mod 5 | i=1,2,3,4 \}
 =\{ 3,5,7,9 \mod 5 \} $
 \par
 Three of $\gamma_i$'s belong to $\Gamma_0(5)$ and
 one of them belongs to $\pm\Gamma_1(5)$.
 Hence
 $(\ref{pole1})
  = \frac{q^{-n}}{pn} \times 2 + O(1)$ and
 $(\ref{pole2})
  = \frac{q^{-n}}{pn} \times 4 + O(1)$.
 \par\noindent
 Case IV. $N=12$ and $p=2$
 \par
 $i=j=1$ and
 $\{ 1+(N/p)ji \mod 2 \}
 =\{ 7 \mod 2 \} $
 \par
 $\gamma_i$ belongs to $\Gamma_0(12)\backslash \pm\Gamma_1(12)$.
 It then follows that
 $(\ref{pole1})
  = \frac{q^{-n}}{pn} + O(1)$ and
 $(\ref{pole2})
  = \frac{q^{-n}}{pn} \times 2 + O(1)$.
 \par\noindent
 Case V. $N=12$ and $p=3$
 \par
 $\{ ji \mod 3 | i=1,2 \}
 =\{ 1,2 \mod 3 \} $
 \par
 $\{ 1+(N/p)ji \mod 3 | i,j=1,2 \}
 =\{ 5,9 \mod 3 \} $
 \par
 Exactly one of $\gamma_i$'s belongs to $\Gamma_0(12)\backslash\pm\Gamma_1(12)$.
 We then have
 $(\ref{pole1})
  = \frac{q^{-n}}{pn} + O(1)$ and
 $(\ref{pole2})
  = \frac{q^{-n}}{pn} \times 2 + O(1)$.
 \par Hence we have the following list of the pole parts.
 \vskip0.3cm
 \begin{center}
 \begin{tabular}{|c||c|c|c|c|}
 \hline
    & $\Gamma_1(N)\infty$ &
    $\Gamma_0(N)\infty\backslash\pm\Gamma_1(N)\infty$ &
    $\frac{1}{(N/p)j}$ ($p\neq 5$) &
    $\frac{1}{(N/p)j}$ ($p=5$) \\
 \hline
 $X_{pn}(t)|_{U_p}$ & $\frac{q^{-n}}{pn}$ & $\times$ & $\frac{q^{-n}}{p^2 n}$
 & $\frac{q^{-n}}{p^2 n}\times 2$ \\
 \hline
 $X_{pn}(t_0)|_{U_p}$ & $\frac{q^{-n}}{pn}$ & $\frac{q^{-n}}{pn}$
 & $\frac{q^{-n}}{p^2 n}\times 2$
 & $\frac{q^{-n}}{p^2 n}\times 4$ \\
 \hline
 $X_n(t)$ & $\frac{q^{-n}}{n}$ & $\times$ & $\times$ & $\times$ \\
 \hline
 $X_n(t_0)$ & $\frac{q^{-n}}{n}$ & $\frac{q^{-n}}{n}$ & $\times$ & $\times$ \\
 \hline
 \end{tabular}
 \end{center}
 \vskip0.3cm
 From the above we see that
 $f(z)=2X_{pn}(t)|_{U_p} - X_{pn}(t_0)|_{U_p}-\frac 1p
 (2X_n(t)-X_n(t_0))$ has no poles. At $\infty$ the Fourier
 expansion of $f(z)$ has no constant term. Thus $f(z)$ is
 identically zero.
 \end{pf}
 \vskip0.2cm
 Write $t_0=\frac 1q+\sum_{m\ge 1} h_m q^m$ and
 $t=\frac 1q+\sum_{m\ge 1} H_m q^m$. Also
 write $X_n(t_0)=\frac 1n q^{-n}+\sum_{m\ge 1} h_{m,n} q^m$ and
 $X_n(t)=\frac 1n q^{-n}+\sum_{m\ge 1} H_{m,n} q^m$.
 \begin{theo}
 Let $N=5,8$ and $p=5$ or $2$ according as $N=5$ or $8$.
 Let $a,b,c,d$ be positive integers such that $ab=cd$ and
 $(a,b)=(c,d)$. Assume that $p^k || (a,b)$. Then
 $$ \psi (a,b) \times (2H_{a,b}-h_{a,b}) =
 \psi (c,d) \times (2H_{c,d}-h_{c,d})$$ where
 $\psi$ is defined by
 $$\psi (a,b)=\begin{cases} 1, & \text{ if } a/p^k \, (\text{ or } b/p^k)
 \equiv \pm 1 \mod N \\ -1, & \text{ otherwise. } \end{cases} $$
 \label{BB}
 \end{theo}
 \begin{pf}
 If $k=0$, then the assertion is immediate from $(*)$. Assume $k\ge 1$. By Lemma \ref{AA},
 $2X_{b}(t)|_{U_p} - X_{b}(t_0)|_{U_p}=\frac 1p
 (2X_{b/p}(t)-X_{b/p}(t_0))$. Comparing the coefficients of
 $q^{a/p}$-terms on both sides, we obtain
 $2H_{a,b}-h_{a,b} =
  \frac 1p (2H_{a/p,b/p}-h_{a/p,b/p})$. By induction
 $2H_{a,b}-h_{a,b} =
  \frac{1}{p^k} (2H_{a/p^k,b/p^k}-h_{a/p^k,b/p^k})$.
 In a similar manner,
 $2H_{c,d}-h_{c,d} =
  \frac{1}{p^k} (2H_{c/p^k,d/p^k}-h_{c/p^k,d/p^k})$.
 Therefore, again by $(*)$, the theorem follows.
 \end{pf}
 \begin{lem} Put $F_{m,n}=\psi (m,n) \times (2H_{m,n}-h_{m,n})$.
 The following two conditions are equivalent:
 \par (i) If $ab=cd$ and $(a,b)=(c,d)$, then $F_{a,b}=F_{c,d}$.
 \par (ii) For any positive integer $s$, there exists an integer $F_m^{(s)}$
 ($0<m\in {\Bbb Z}$) such that $F_m^{(1)}=2H_m-h_m$
 and $F_{m,n}=\sum_{s|(m,n)} \frac 1s F_{mn/s^2}^{(s)}$.
 (Convention: when a summation runs over the divisors of a number,
 these are considered to be positive.)
 \label{CC}
 \end{lem}
 \begin{pf}
 (ii) $\Rightarrow$ (i) clear.
 \par\noindent
 (i) $\Rightarrow$ (ii)
 Let $F$ be a numerical function. If we define $f$ by
 $f(n)=\sum_{d|n} \mu (d) F(n/d)$ ($n=1,2,3,\cdots$), then
 by the M\"obius inversion formula
 $F(n)=\sum_{d|n} f(n/d)$ for all $n\in {\Bbb Z}$. The converse
 holds, too. For a fixed $mn$, $F_{m,n}$ can be viewed as a
 function of $(m,n)$. Define $F_k^{(s)}=s\cdot
 \sum_{d|s} \mu (d) F_{m,n}$ for $s,k\ge 1$ where we choose $m,n$ so that
 $mn=s^2 k$ and $(m,n)=s/d$ (e.g. $m=s/d$ and $n=(s/d)\cdot d^2
 k$). Then by the condition (i), $F_k^{(s)}$ is well-defined and
 $F_m^{(1)}=1\cdot \mu (1) \cdot F_{m,1}=2H_m-h_m$.
 Here $\psi (m,1)=1$ for $N=5,8$ (as described in Theorem \ref{BB}).
 If $N=10$ or $12$, then $\psi (m,1)$ will be defined to be $1$
 (see \S 3 and \S 4).
 Define a map $F:
 {\Bbb N} \to {\Bbb C}$ by
 $F(x)=\begin{cases} F_{a,b}, & \text{ if } ab=mn \text{ and }
 x=(a,b) \\ 0, & \text{ otherwise. } \end{cases}$
 \par\noindent Then
 $\sum_{s|(m,n)} \frac 1s
 F_{mn/s^2}^{(s)}=\sum_{s|(m,n)}\frac 1s \cdot s \cdot \sum_{d|s}
 \mu (d) F_{\frac sd, \frac sd \cdot d^2 \cdot \frac{mn}{s^2}}=
 \sum_{s|(m,n)}\sum_{d|s}\mu (d) F(s/d)=F((m,n))=F_{m,n}$.
 \end{pf}
 \begin{rmk} For a negative integer $m$ and a positive integer $n$, we define
 \begin{align*}
  & F_{m,n}=\psi (m,n)=H_{m,n}=h_{m,n}=
   \begin{cases} 1, & \text{ if } m=-n \\
                 0, & \text{ if } m\neq -n
   \end{cases} \\
   \intertext{ and }
  & F_m^{(s)}=\begin{cases} 1, & \text{ if } m=-1 \\
                 0, & \text{ if } m < -1.
   \end{cases}
 \end{align*}
 Then Lemma \ref{CC}-(ii) can be extended to all $m\in {\Bbb Z}$.
 \end{rmk}

 \begin{prop}
 Let
 $p=e^{2\pi i y}$ and $q=e^{2\pi i z}$ with $y,z\in {\frak H}$.
 Then $F_m^{(s)}$ satisfy the following product identity:
 $$ p^{-1}\prod_{n>0\atop m\in {\Bbb Z}} \exp\left(-\sum_{s\ge 1}
  \left(\frac{\psi(sm,sn) F_{mn}^{(s)}+
  h_{mn}^{(s)}}{2s}\right)
    q^{sm} p^{sn}\right)= t(p)-t(q). $$
 \label{DD}
 \end{prop}
 \begin{pf}
 Note that $X_n(t)$ can be viewed as the coefficient of $p^n$-term
 in $-\log p - \log (t(p)-t(q))$ (see \cite{Norton}). Thus
 $\log p^{-1}-\sum_{n > 0} X_n(t) p^n = \log (t(p)-t(q))$.
 Taking exponential on both sides, we get
 \begin{equation}
 p^{-1} \exp (-\sum_{n> 0} X_n(t) p^n)=t(p)-t(q).
 \label{product1}
 \end{equation}
 On the other hand,
 \begin{align*}
 \sum_{n>0} X_n(t) p^n
 & =\sum_{n>0}\sum_{m \in {\Bbb Z}} H_{m,n} q^m p^n
   =\sum_{n>0}\sum_{m \in {\Bbb Z}} \frac{\psi
   (m,n)F_{m,n}+h_{m,n}}{2} q^m p^n \\
 & =\sum_{n>0}\sum_{m \in {\Bbb Z}}\sum_{s|(m,n)}
   \left(\frac{\psi(m,n) F_{mn/s^2}^{(s)}+
  h_{mn/s^2}^{(s)}}{2s}\right) q^m p^n \\
 & =\sum_{n>0}\sum_{m \in {\Bbb Z}}\sum_{s=1}^\infty
   \left(\frac{\psi(sm,sn) F_{mn}^{(s)}+
  h_{mn}^{(s)}}{2s}\right) q^{sm} p^{sn}.
 \end{align*}
 By plugging the above into (\ref{product1}) we obtain the desired
 product identity.
 \end{pf}
 \begin{theo} (Characterization of super-replicable functions)
 Let $t=\frac 1q+\sum_{m\ge 1} H_m q^m$ be a normalized $q$-series . The
 followings are equivalent.
 \par (a) $t$ is super-replicable with respect to $t_0=\frac
 1q+\sum_{m\ge 1} h_m q^m$ and a character $\psi$.
 \par (b) for any $s\ge 1$ there exists a normalized $q$-series $F^{(s)}=\frac
 1q+\sum_{m\ge 1} F_m^{(s)} q^m$ such that
 $ F^{(1)}=2t-t_0$ and $$(2X_n(t)-X_n(t_0))_{\psi}
  =n^{-1}\sum_{ad=n \atop 0\le b < d}
 F^{(a)}\left(\frac{az+b}{d}\right)$$
 where the notation $(2X_n(t)-X_n(t_0))_{\psi}$ is used to mean
 a $q$-series $\sum_{m\in {\Bbb Z}} \psi (m,n) (2H_{m,n}-h_{m,n})
 q^m$.
 \par (c) for every $s\ge 1$ there exists a normalized $q$-series $F^{(s)}=\frac
 1q+\sum_{m\ge 1} F_m^{(s)} q^m$ such that
 $$ p^{-1}\prod_{n>0\atop m\in {\Bbb Z}} \exp\left(-\sum_{s\ge 1}
  \left(\frac{\psi(sm,sn) F_{mn}^{(s)}+
  h_{mn}^{(s)}}{2s}\right)
    q^{sm} p^{sn}\right)= t(p)-t(q) $$ where $p=e^{2\pi i y}$ and $q=e^{2\pi i
    z}$ with $y,z\in {\frak H}$.
 \label{product2}
 \end{theo}
 \begin{pf}
 The equivalence of (a) and (b) follows from Lemma \ref{CC}.
 And Proposition \ref{DD} yields that (a) implies (c).
 \par\noindent (c) $\Rightarrow$ (a)
 \par\noindent The identity in (c) can be rewritten as
 $$ p^{-1}\prod_{n>0\atop m\in {\Bbb Z}} \exp\left(-\sum_{s|(m,n)}
  \left(\frac{\psi(m,n) F_{mn/s^2}^{(s)}+
  h_{mn/s^2}^{(s)}}{2s}\right)
    q^{m} p^{n}\right)= t(p)-t(q). $$
 On the other hand,
 $t(p)-t(q)=p^{-1} \exp(-\sum_{n>0} X_n(t) p^n)
           =p^{-1} \exp(-\sum_{n>0 \atop m\in{\Bbb Z}} H_{m,n}q^m
           p^n)$.
 We then have $H_{m,n}=-\sum_{s|(m,n)}
  \left(\frac{\psi(m,n) F_{mn/s^2}^{(s)}+
  h_{mn/s^2}^{(s)}}{2s}\right)$, which implies
  $\psi (m,n) (2H_{m,n}-h_{m,n})
  = \sum_{s|(m,n)}F_{mn/s^2}^{(s)}$ because $h$ is replicable.
  The last identity then yields (a).
 \end{pf}

 \section{The case $N=12$}
 \par
 \begin{lem}
 For odd $r$,
 \par (i) $(X_r(t)|_{U_{2^k}})|_{\sm 1&0\\ 6&1 \esm}
 =\begin{cases} \frac{1}{2r} q_2^{-r}+O(1), & \text{ if } k=1 \\
                O(1), & \text{ otherwise }
  \end{cases} $
 \par (ii) $(X_r(t)|_{U_{2^k}})|_{\sm 1&0\\ 3&1 \esm}
 =\begin{cases} O(1), & \text{ if } k=1 \\
                \frac{1}{4r} q_4^{-r}+O(1), & \text{ if } k=2 \\
                \frac{1}{8r} i^{-r}\cdot q_2^{-r}+O(1), & \text{ if } k=3 \\
                \frac{(-1)^k i^{-r}}{2^k r} q^{-2^{k-4}r}+O(1), & \text{ if } k\ge
                4.
  \end{cases} $
 \label{EE}
 \end{lem}
 \begin{pf}
 (i) First, note that ${X_r(t)}|_{\sm 1&0\\6&1 \esm}$ is holomorphic at
$\infty$ because ${X_r(t)}$ has poles only at the cusps
$\Gamma_1(12)\infty$. Now for $k\ge 1$,
\begin{align*}
\left({X_r(t)}|_{U_{2^k}}\right)|_{\sm 1&0\\6&1 \esm }
&=\left(\left({X_r(t)}|_{U_{2^{k-1}}}\right)|_{U_2}\right)|_{\sm
1&0\\6&1 \esm }
\\
&=\frac12\left({X_r(t)}|_{U_{2^{k-1}}}\right)|_{\sm 1&0\\0&2 \esm
\sm 1&0\\6&1 \esm}
+\frac12\left({X_r(t)}|_{U_{2^{k-1}}}\right)|_{\sm 1&1\\0&2 \esm
\sm 1&0\\6&1 \esm}
\\
&=\frac12\left({X_r(t)}|_{U_{2^{k-1}}}\right)|_{\sm 1&0\\12&1 \esm
\sm 1&0\\0&2 \esm}
+\frac12\left({X_r(t)}|_{U_{2^{k-1}}}\right)|_{\sm 7&4\\12&7 \esm
\sm 1&-1\\0&2 \esm}
\\
&=\frac12\left({X_r(t)}|_{U_{2^{k-1}}}\right)\left(\frac z2\right)
+\frac12\left(\left({X_r(t)}|_{\sm 7&4\\12&7 \esm}\right)|_
{U_{2^{k-1}}}\right)|_{\sm 1&-1\\0&2 \esm}
\text{ \q by Lemma \ref{Hecke} } \\
&=\begin{cases}
\frac{1}{2r} {q_2}^{-r}+O(1), & \text{ if $k=1$} \\
O(1), & \text{ otherwise. }
\end{cases}
\end{align*}
\par\noindent
(ii) We observe that ${X_r(t)}|_{\sm 1&0\\3&1 \esm}\in O(1)$. And
for $k\ge 1$,
\begin{align}
\left({X_r(t)}|_{U_{2^k}}\right)|_{\sm 1&0\\3&1 \esm }
&=\left(\left({X_r(t)}|_{U_{2^{k-1}}}\right)|_{U_2}\right)|_{\sm
1&0\\3&1 \esm }
\label{HH} \\
& =\frac12\left({X_r(t)}|_{U_{2^{k-1}}}\right)|_{\sm 1&0\\0&2 \esm
\sm 1&0\\3&1 \esm}
+\frac12\left({X_r(t)}|_{U_{2^{k-1}}}\right)|_{\sm 1&1\\0&2 \esm
\sm 1&0\\3&1 \esm}
\notag \\
&=\frac12\left({X_r(t)}|_{U_{2^{k-1}}}\right)|_{\sm 1&0\\6&1 \esm
\sm 1&0\\0&2 \esm}
+\frac12\left({X_r(t)}|_{U_{2^{k-1}}}\right)|_{\sm 2&1\\3&2 \esm
\sm 2&0\\0&1 \esm} \notag
\end{align}
has a holomorphic Fourier expansion if $k=1$. Thus we suppose
$k\ge 2$. We then derive that
\begin{align*}
\left({X_r(t)}|_{U_{2^{k-1}}}\right)|_{\sm 2&1\\3&2 \esm }
&=\left(\left({X_r(t)}|_{U_{2^{k-2}}}\right)|_{U_2}\right)|_{\sm 2&1\\3&2 \esm } \\
&=\frac12\left({X_r(t)}|_{U_{2^{k-2}}}\right)|_{\sm 1&0\\0&2 \esm
\sm 2&1\\3&2 \esm}
+\frac12\left({X_r(t)}|_{U_{2^{k-2}}}\right)|_{\sm 1&1\\0&2 \esm
\sm 2&1\\3&2 \esm}
\\
&=\frac12\left({X_r(t)}|_{U_{2^{k-2}}}\right)|_{\sm 1&0\\3&1 \esm
\sm 2&1\\0&1 \esm} +\frac12\left({X_r(t)}|_{U_{2^{k-2}}}\right)|_
{\sm 11&-1\\12&-1 \esm \sm 1&0\\6&1 \esm\sm 1&1\\0&2 \esm}
\\
&=\frac12\left({X_r(t)}|_{U_{2^{k-2}}}\right)|_{\sm 1&0\\3&1 \esm
\sm 2&1\\0&1 \esm}
+\frac12\left({X_r(t)}|_{U_{2^{k-2}}}\right)|_{\sm 1&0\\6&1 \esm
\sm 1&1\\0&2 \esm} \, .
\end{align*}
If we substitute the above for (\ref{HH}), we get that for $k\ge
2$,
\begin{align}
 & \q  \left({X_r(t)}|_{U_{2^{k}}}\right)|_{\sm 1&0\\3&1 \esm } \label{IIII}\\
 & =\frac12\left({X_r(t)}|_{U_{2^{k-1}}}\right)|_{\sm 1&0\\6&1 \esm
   \sm 1&0\\0&2 \esm }
  +\frac14\left({X_r(t)}|_{U_{2^{k-2}}}\right)|_{\sm 1&0\\6&1 \esm \sm 1&1\\0&2 \esm \sm 2&0\\0&1 \esm }
  +\frac14\left({X_r(t)}|_{U_{2^{k-2}}}\right)|_{\sm 1&0\\3&1 \esm \sm 2&1\\0&1 \esm \sm 2&0\\0&1 \esm }
   \notag \\
 & =\frac12\left({X_r(t)}|_{U_{2^{k-1}}}\right)|_{\sm 1&0\\6&1 \esm
   \sm 1&0\\0&2 \esm }
 +\frac14\left({X_r(t)}|_{U_{2^{k-2}}}\right)|_{\sm 1&0\\6&1 \esm \sm 2&1\\0&2 \esm }
 +\frac14\left({X_r(t)}|_{U_{2^{k-2}}}\right)|_{\sm 1&0\\3&1 \esm \sm 4&1\\0&1 \esm
 } \, .
\notag
\end{align}
When $k=2$,
\begin{align*}
& \q  \left({X_r(t)}|_{U_{4}}\right)|_{\sm 1&0\\3&1 \esm }
=\frac12\left({X_r(t)}|_{U_{2}}\right)|_{\sm 1&0\\6&1
\esm}\left(\frac z2\right)
 +\frac14 {X_r(t)}|_{\sm 1&0\\6&1 \esm }\left(z+\frac12\right)
 +\frac14 {X_r(t)}|_{\sm 1&0\\3&1 \esm }\left(4z+1\right)  \\
& = \frac{1}{4r} {q_4}^{-r}+O(1) \text{ \q\q by (i). }
\end{align*}
If $k=3$,
\begin{align*}
& \q  \left({X_r(t)}|_{U_{8}}\right)|_{\sm 1&0\\3&1 \esm } \\
&=\frac12\left({X_r(t)}|_{U_{4}}\right)|_{\sm 1&0\\6&1
\esm}\left(\frac z2\right)
 +\frac14\left({X_r(t)}|_{U_{2}}\right)|_{\sm 1&0\\6&1 \esm } \left(z+\frac12\right)
 +\frac14\left({X_r(t)}|_{U_{2}}\right)|_{\sm 1&0\\3&1 \esm }(4z+1) \\
& = \frac{1}{8r} e^{-\pi i\left(z+\frac12\right)r} + O(1) =
\frac{i^{-r}}{8r}{q_2}^{-r}+O(1) \text{ \q\q by (i) and the case
$k=1$ in (ii). }
\end{align*}
For $k\ge 4$, we will show by induction on $k$ that
\begin{equation}
\left({X_r(t)}|_{U_{2^{k}}}\right)|_{\sm 1&0\\3&1 \esm }
   = (-1)^{k}\cdot \frac{i^{-r}}{2^k r}\cdot q^{-2^{k-4} r}+O(1).
   \label{JJJJ}
\end{equation}
First we note that by (i) and (\ref{IIII})
$$\left({X_r(t)}|_{U_{2^{k}}}\right)|_{\sm 1&0\\3&1 \esm }
  =\frac14\left({X_r(t)}|_{U_{2^{k-2}}}\right)|_{\sm 1&0\\3&1 \esm }(4z+1)+O(1). $$
If $k=4$,
\begin{align*}
\left({X_r(t)}|_{U_{2^{4}}}\right)|_{\sm 1&0\\3&1 \esm }
  &=\frac14\left({X_r(t)}|_{U_{2^{4-2}}}\right)|_{\sm 1&0\\3&1 \esm }(4z+1)+O(1)\\
&=\frac{1}{16 r}e^{-\frac{\pi i}{2}(4z+1)r}+O(1)
  \text{ \q \q by the case $k=2$} \\
&= \frac{i^{-r}}{2^4r}\cdot q^{-r}+O(1).
\end{align*}
Thus when $k=4$, (\ref{JJJJ}) holds. Meanwhile, if $k=5$ then we
get that
\begin{align*}
\left({X_r(t)}|_{U_{2^{5}}}\right)|_{\sm 1&0\\3&1 \esm }
  &=\frac14\left({X_r(t)}|_{U_{2^{5-2}}}\right)|_{\sm 1&0\\3&1 \esm }(4z+1)+O(1)\\
&=\frac{i^{-r}}{32r}e^{-\pi i(4z+1)r}+O(1)
  \text{ \q \q by the case $k=3$} \\
&= (-1)\cdot \frac{i^{-r}}{2^5 r}\cdot q^{-2r}+O(1).
\end{align*}
Therefore in this case (\ref{JJJJ}) is also valid. Now for $k\ge
6$,
\begin{align*}
\left({X_r(t)}|_{U_{2^{k}}}\right)|_{\sm 1&0\\3&1 \esm }
  &=\frac14\left({X_r(t)}|_{U_{2^{k-2}}}\right)|_{\sm 1&0\\3&1 \esm }(4z+1)+O(1)\\
&=\frac{1}{4}\cdot (-1)^{k-2}\cdot \frac{i^{-r}}{2^{k-2}r}\cdot
 e^{-2\pi i(4z+1)\cdot 2^{k-2-4}r}+O(1) \\
& \text{ \q \q \, by induction hypothesis for $k-2$} \\
&= (-1)^{k}\cdot \frac{i^{-r}}{2^kr}\cdot q^{-2^{k-4}r}+O(1).
\end{align*}
This proves the lemma.
\end{pf}

 For a modular function $f=\sum_{n\in {\Bbb Z}} a_n q^n$ and
 $\chi (n)=\left(\frac{-1}{n}\right)$ (the Jacobi symbol), the twist of $f$ by $\chi$
 is defined by
 $f_{\chi}=\sum_{n\in {\Bbb Z}} a_n \chi (n) q^n =\frac{1}{\sqrt{-4}}
  (f(z+1/4)-f(z+3/4))$ (see \cite{Koblitz}, pp. 127-128).
 \begin{lem}
 For odd $r$,
 $$ (X_r(t)|_{U_{2^k}})_\chi (z)=(-1)^k\cdot i^{-r-1}\cdot \frac
 12 (X_{2^k r}(t)(z)-X_{2^k r}(t)(z+1/2)).$$
 \label{FF}
 \end{lem}
 \begin{pf} Denote the right hand side by $g(z)$. We observe that
 $g(z)=(-1)^k\cdot i^{-r-1}\cdot (X_{2^k r}(t)(z)-X_{2^k
 r}(t)|_{U_2}(2z))$. Thus $g(z)\in K(X_1(24))$.
 In the proof of \cite{Super}, Lemma 15 and 16, it is shown that for
 $f\in K(X_1(12))$ and $\sm a&b\\ c&d \esm\in \Gamma_0(12)$ such that
 $(f-c)|_{\sm a&b\\ c&d \esm}=\left(\frac 3a \right)\cdot
 (f-c)$ with some constant $c$ \, $f_\chi$ belongs to $K(X_1(24))$.
 And it is also proved there that for $\sm a&b\\ c&d \esm\in \Gamma_0(24)$,
 ${f_\chi}|_{\sm a&b\\ c&d \esm}=(-1)^{c/24}\cdot
 \left(\frac{3}{a+\frac c4}\right)\cdot f_\chi$.
 We take $X_r(t)|_{U_{2^k}}$ in place of $f$.
 Then we derive that for $\sm a&b\\ c&d \esm\in
 \Gamma_0(12)$,
 \begin{align*}
 & (X_r(t)|_{U_{2^k}})|_{\sm a&b\\ c&d \esm}
  = (X_r(t)|_{\sm a&b\\ c&d \esm})|_{U_{2^k}}
  = \frac 12 \{\psi_{12}
 (d)(2X_r(t)-X_r(t_0))+X_r(t_0)\}|_{U_{2^k}}+c_0 \\
 & \text{ with } c_0=
   \begin{cases} 0, & \text{ if } d\equiv \pm 1 \mod 12 \\
                 X_r(t)(a/c), & \text{ otherwise }
   \end{cases} \\
 & \text{ by Lemma \ref{Hecke} and \ref{Action} }
   \\
 & = \psi_{12} (d) X_r(t)|_{U_{2^k}}+c_0
    \text{ since } X_r(t_0)|_{U_{2^k}}\equiv 0 \text{ by \cite{Koike},
    Theorem 3.1 Case V} \\
 & = \left( \frac 3a \right) X_r(t)|_{U_{2^k}}+c_0.
 \end{align*}
 Thus $(X_r(t)|_{U_{2^k}}-c_0/2)|_{\sm a&b\\ c&d \esm}
 =\left(\frac 3a \right)\cdot
 (X_r(t)|_{U_{2^k}}-c_0/2)$, which leads to the fact that
 $(X_r(t)|_{U_{2^k}})_\chi \in K(X_1(24))$ and
 for $\sm a&b\\ c&d \esm\in \Gamma_0(24)$,
 \begin{equation}
 {(X_r(t)|_{U_{2^k}})_\chi}|_{\sm a&b\\ c&d \esm}=(-1)^{c/24}\cdot
 \left(\frac{3}{a+\frac c4}\right)\cdot (X_r(t)|_{U_{2^k}})_\chi.
 \label{twist}
 \end{equation}
 Here we note that both $(X_r(t)|_{U_{2^k}})_\chi$ and $g$ sit in $K(X_1(24))$.
 As for the equality
 $(X_r(t)|_{U_{2^k}})_\chi=g$, it suffices to show that $(X_r(t)|_{U_{2^k}})_\chi-g$
 has no poles in ${\frak H}^*$. By the same argument as in
 \cite{Super}, Lemma 18, $(X_r(t)|_{U_{2^k}})_\chi-g$ can have poles
 only at the cusps which are equivalent to $1/12$ and $5/12$ under
 $\Gamma_1(24)$.
 \par\noindent
 At $1/12$,
 \begin{align*}
 \left(X_r(t)|_{U_{2^k}}\right)_\chi|_{\sm 1&0\\12&1 \esm }
 &=\frac{1}{\sqrt{-4}} \left(X_r(t)|_{U_{2^k}\sm 4&1\\0&4 \esm\sm
 1&0\\12&1 \esm}
     -X_r(t)|_{U_{2^k}\sm 4&3\\0&4 \esm\sm 1&0\\12&1 \esm}\right)
  \\
 &=\frac{1}{\sqrt{-4}} \left(X_r(t)|_{U_{2^k}\sm 1&0\\3&1 \esm\sm
 16&1\\0&1 \esm}
     -X_r(t)|_{U_{2^k}\sm 5&-1\\6&-1 \esm\sm 8&1\\0&2 \esm}\right)  \\
 &=\frac{1}{\sqrt{-4}} \left(X_r(t)|_{U_{2^k}\sm 1&0\\3&1 \esm\sm
 16&1\\0&1 \esm}
     -X_r(t)|_{U_{2^k}\sm 11&-1\\12&-1 \esm\sm 1&0\\6&1 \esm
                 \sm 8&1\\0&2 \esm}\right)
  \\
 &=\frac{1}{\sqrt{-4}} \left(X_r(t)|_{U_{2^k}\sm 1&0\\3&1 \esm}(16z+1)
     -X_r(t)|_{U_{2^k}\sm 1&0\\6&1
     \esm}\left(4z+\frac12\right)\right),
 \end{align*}
 which is, by Lemma \ref{EE}, equal to
 \begin{align*}
 & \begin{cases}
  \frac{1}{\sqrt{-4}}\left(-\frac{1}{2r} e^{-\pi i r(4z+1/2)} +
  O(1)\right), & \text{ if } k=1 \\
  \frac{1}{\sqrt{-4}}\left(\frac{1}{4r} e^{-\frac{\pi i}{2} r(16z+1)} +
  O(1)\right), & \text{ if } k=2 \\
  \frac{1}{\sqrt{-4}}\left(\frac{1}{8r}\cdot i^{-r} e^{-\pi i r(16z+1)} +
  O(1)\right), & \text{ if } k=3 \\
  \frac{1}{\sqrt{-4}}\left(\frac{(-1)^k\cdot i^{-r}}{2^k\cdot r}
  e^{-2\pi i r\cdot 2^{k-4}(16z+1)} +
  O(1)\right), & \text{ if } k\ge 4
  \end{cases} \\
  &=(-1)^k\cdot i^{-r-1}\cdot \frac{1}{2^{k+1} r}\cdot q^{-r\cdot
  2^k} + O(1) \text{ for } k\ge 1.
 \end{align*}
 On the other hand,
 \begin{align*}
 g|_{\sm 1&0\\12&1 \esm} & = (-1)^{k}\cdot i^{-r-1}\cdot \frac12\cdot \left(
 X_{2^k r}(t)|_{\sm 1&0\\12&1 \esm}
  -X_{2^k r}(t)|_{\sm 2&1\\0&2 \esm\sm 1&0\\12&1 \esm}\right) \\
 & = (-1)^{k}\cdot i^{-r-1}\cdot \frac12\cdot \left( X_{2^k r}(t)
  -X_{2^k r}(t)|_{\sm 7&4\\12&7 \esm\sm 2&-1\\0&2 \esm}\right) \\
 & = (-1)^{k}\cdot i^{-r-1}\cdot \frac{1}{2^{k+1}r}\cdot q^{-2^k r}+O(1).
 \end{align*}
 At $5/12$,
 \begin{align*}
 \left(X_r(t)|_{U_{2^k}}\right)_\chi|_{\sm 5&2\\12&5 \esm}
 &= \left(X_r(t)|_{U_{2^k}}\right)_\chi|_
  {\sm -19&2\\-48&5 \esm\sm 1&0\\12&1 \esm} \\
 &= (-1)^\frac{-48}{24}\cdot
 \left(\frac{3}{-19+\frac{-48}{4}}\right)
  \cdot\left(X_r(t)|_{U_{2^k}}\right)_\chi|_{\sm 1&0\\12&1 \esm}
  \text{\q by (\ref{twist})} \\
 &= (-1)\cdot \left(X_r(t)|_{U_{2^k}}\right)_\chi|_{\sm 1&0\\12&1
 \esm }
 = (-1)^{k+1}\cdot i^{-r-1}\frac{1}{2^{k+1}r}\cdot q^{-r\cdot 2^k}+O(1).
 \end{align*}
 And
 \begin{align*}
 g|_{\sm 5&2\\12&5 \esm} & = (-1)^{k}\cdot i^{-r-1}\cdot \frac12\cdot \left(
 X_{2^k r}(t)|_{\sm 5&2\\12&5 \esm}
  -X_{2^k r}(t)|_{\sm 2&1\\0&2 \esm\sm 5&2\\12&5 \esm}\right) \\
 & = (-1)^{k}\cdot i^{-r-1}\cdot \frac12\cdot \left( X_{2^k}(t)|_{\sm
 5&2\\12&5 \esm}
  -X_{2^k}(t)|_{\sm 11&-1\\12&-1 \esm\sm 2&1\\0&2 \esm}\right) \\
 & = (-1)^{k+1}\cdot i^{-r-1}\cdot \frac{1}{2}\cdot X_{2^k r}(t)|_{\sm 2&1\\0&2
 \esm}+O(1) = (-1)^{k+1}\cdot i^{-r-1}\cdot \frac{1}{2^{k+1}r}\cdot q^{-2^k r}+O(1).
 \end{align*}
 This implies that
 $\left(X_r(t)|_{U_{2^k}}\right)_\chi-g$ has
 no poles, as desired.
 Now the lemma is proved.
 \end{pf}
 Let $a$ and $b$ be positive integers.
 We define $\psi (a,b)$ as follows.
 First we assume that $(6,(a,b))=1$.
 \par\noindent
 Case I: $a$ odd and $b$ odd.
 $$\psi (a,b)= \begin{cases} 1, & \text{ if } a\equiv \pm 1
 \mod 12 \text{ \, or \, } b\equiv \pm 1 \mod 12 \\
 -1, & \text{ otherwise.} \end{cases}$$
 \par\noindent
 Case II: $a$ even and $b$ odd. Write $a=2^k n$ with
 $n$ odd.
 $$\psi (a,b)= \begin{cases} \psi_{12} (b), & \text{ if } (6,b)=1 \\
 (-1)^k\cdot i^{b+1} \cdot \chi (n) \cdot \psi_{12} (n), & \text{ otherwise}
 \end{cases}$$
 where $\chi (n)=\left( \frac{-1}{n} \right)$.
 \par\noindent
 Case III: $a$ odd and $b$ even.
 \par Put $\psi (a,b)= \psi (b,a)$.
 \par\noindent
 In general, let $2^\alpha 3^\beta || (a,b)$ and define
 $\psi (a,b) = \psi (\frac{a}{2^\alpha 3^\beta},\frac{b}{2^\alpha
 3^\beta})$.
 \begin{theo}
 Let $N=12$. If $ab=cd$ and $(a,b)=(c,d)$, then
 $$ \psi (a,b) \times (2H_{a,b}-h_{a,b}) =
 \psi (c,d) \times (2H_{c,d}-h_{c,d}).$$
 \label{GG}
 \end{theo}
 \begin{pf}
 Let $2^\alpha 3^\beta || (a,b)=(c,d)$.
 By Lemma \ref{AA},
 \begin{align*}
  & 2H_{a,b}-h_{a,b}=\frac{1}{2^\alpha 3^\beta}
    (2H_{\frac{a}{2^\alpha 3^\beta},\frac{b}{2^\alpha 3^\beta}}
     -h_{\frac{a}{2^\alpha 3^\beta},\frac{b}{2^\alpha 3^\beta}})
     \\
  \intertext{ and }
  & 2H_{c,d}-h_{c,d}=\frac{1}{2^\alpha 3^\beta}
    (2H_{\frac{c}{2^\alpha 3^\beta},\frac{d}{2^\alpha 3^\beta}}
     -h_{\frac{c}{2^\alpha 3^\beta},\frac{d}{2^\alpha 3^\beta}}).
 \end{align*}
 Thus we may assume that $(6,(a,b))=(6,(c,d))=1$.
 \par\noindent
 Case I: $a$ odd and $b$ odd
 \par
 In this case $(a,12)=1$ or $(b,12)=1$. And $(c,12)=1$ or
 $(d,12)=1$. Therefore the assertion follows from $(*)$.
 \par\noindent
 Case II: $a$ even or $b$ even
 \par
 Since $\psi (a,b)=\psi (b,a)$, $H_{a,b}=H_{b,a}$ and
 $h_{a,b}=h_{b,a}$, we may assume that $a$ even and $b$ odd
 (resp. $c$ even and $d$ odd).
 Write $a=2^k\cdot a'$ with $a'$ odd (resp.
 $c=2^k\cdot c'$ with $c'$ odd) and $k\ge 1$.
 By Lemma \ref{FF} it holds that
 $$ (X_b(t)|_{U_{2^k}})_\chi (z)=(-1)^k\cdot i^{-b-1}\cdot \frac
 12 (X_{2^k b}(t)(z)-X_{2^k b}(t)(z+1/2)).$$
 Comparing the coefficients of $q^{a'}$-terms on both sides, we
 obtain that $\chi (a') H_{2^k a',b}=(-1)^k\cdot
 i^{-b-1}\cdot H_{a',2^k b}$.
 Since $h_{2^k a',b}=h_{a',2^k b}=0$ by \cite{Koike}, Theorem 3.1 Case
 V, we end up with
 $\chi (a')\cdot (-1)^k\cdot i^{b+1}\cdot
  (2H_{a,b}-h_{a,b})=2H_{a',2^k b}-h_{a',2^k b}$.
 Similarly,
 $\chi (c')\cdot (-1)^k\cdot i^{d+1}\cdot
  (2H_{c,d}-h_{c,d})=2H_{c',2^k d}-h_{c',2^k d}$.
 \par\noindent
 If $(b,12)=(d,12)=1$, then by $(*)$
 $$ \psi_{12} (b) \times (2H_{a,b}-h_{a,b}) =
 \psi_{12} (d) \times (2H_{c,d}-h_{c,d}).$$
 Here we note by definition that $\psi_{12} (b)=\psi (a,b)$ and
 $\psi_{12} (d)=\psi (c,d)$.
 \par\noindent
 If $(b,12)=1 $ and $(d,12)\neq 1$, then $(c',12)=1$. Again
 by $(*)$
 \begin{align*}
  \psi_{12} (b) \times (2H_{a,b}-h_{a,b})
  &= \psi_{12} (c') \times (2H_{c',2^k d}-h_{c',2^k d}) \\
  &= (-1)^k\cdot i^{d+1}\chi (c') \psi_{12} (c') \cdot
     (2H_{c,d}-h_{c,d}).
 \end{align*}
 Now the assertion follows from the definition of $\psi$.
 In the remaining cases (i.e. $(b,12)\neq 1$ and
 $(d,12)=1$ or $(b,12)\neq 1$ and $(d,12)\neq 1$)
 the same argument leads to the assertion.
 \end{pf}
 \section{The case $N=10$}
 \par
 Throughout this section we make use of the following notations.
 \par
 $t$: Hauptmodul of $\Gamma_1(10)$, \,
 $t^{(2)}$: Hauptmodul of $\Gamma_1(5)$
 \par
 $t_0$: Hauptmodul of $\Gamma_0(10)$, \,
 $t_0^{(2)}$: Hauptmodul of $\Gamma_0(5)$
 \begin{lem}
 $2X_{2n}(t)|_{U_2} - X_{2n}(t_0)|_{U_2}=\frac 14
 (2X_n(t^{(2)})-X_n(t_0^{(2)}))
 +\frac 14 (2X_n(t)-X_n(t_0))$.
 \label{II}
 \end{lem}
 \begin{pf}
 First we notice that all functions in the assertion are invariant
 under $\Gamma_1(10)$.
 Now we investigate their possible poles.
 \par\noindent
 $X_n(t^{(2)})$ (resp. $X_n(t_0^{(2)})$) can have poles only at $\Gamma_1(5)\infty$
 (resp. $\Gamma_0(5)\infty$). And
 $X_n(t)$ (resp. $X_n(t_0)$) can have poles only at $\Gamma_1(10)\infty$
 (resp. $\Gamma_0(10)\infty$). Furthermore
 $X_{2n}(t)|_{U_2}$ (resp. $X_{2n}(t_0)|_{U_2}$) can have poles only at $\sm 1&i \\ 0&2
 \esm^{-1}\Gamma_1(10)\infty$ (resp. $\sm 1&i \\ 0&2
 \esm^{-1}\Gamma_0(10)\infty$) for $i=0,1$.
 \par\noindent
 Then, up to $\Gamma_1(10)$-equivalence,
 \par\noindent
 cusps in $\Gamma_1(5)\infty$ are equivalent to $1/5$ or $\infty$,
 \par\noindent
 cusps in $\Gamma_0(5)\infty$ are equivalent to $1/5$ or $2/5$ or $\infty$ or
 $3/10$,
 \par\noindent
 cusps in $\Gamma_1(10)\infty$ are equivalent to $\infty$,
 \par\noindent
 cusps in $\Gamma_0(10)\infty$ are equivalent to $\infty$ or
 $3/10$,
 \par\noindent
 cusps in $\sm 1&i \\ 0&2 \esm^{-1}\Gamma_1(10)\infty$ are equivalent to
 $1/5$ or $\infty$, and
 \par\noindent
 cusps in $\sm 1&i \\ 0&2 \esm^{-1}\Gamma_0(10)\infty$ are equivalent to
 $1/5$ or $2/5$ or $\infty$ or $3/10$ (see the proof of Lemma
 \ref{AA}).
 \par
 By investigating the pole parts at the cusps mentioned above we have the
 following table:
 \vskip0.3cm
 \begin{center}
 \begin{tabular}{|c||c|c|c|c|}
 \hline
  & $\infty$ & $3/10$ & $2/5$ & $1/5$ \\
 \hline
 $X_{2n}(t)|_{U_2}$ & $q^{-n}/(2n)$ & $\times$ & $\times$ &
 $q^{-n}/(4n)$ \\
 \hline
 $X_{2n}(t_0)|_{U_2}$ & $q^{-n}/(2n)$ & $q^{-n}/(2n)$ & $q^{-n}/(4n)$ &
 $q^{-n}/(4n)$ \\
 \hline
 $\frac 14 X_n(t^{(2)})$ & $q^{-n}/(4n)$ & $\times$ & $\times$ &
 $q^{-n}/(4n)$ \\
 \hline
 $\frac 14 X_n(t_0^{(2)})$ & $q^{-n}/(4n)$ & $q^{-n}/(4n)$ & $q^{-n}/(4n)$ &
 $q^{-n}/(4n)$ \\
 \hline
 $\frac 14 X_n(t)$ & $q^{-n}/(4n)$ & $\times$ & $\times$ &
 $\times$ \\
 \hline
 $\frac 14 X_n(t_0)$ & $q^{-n}/(4n)$ & $q^{-n}/(4n)$ & $\times$ &
 $\times$ \\
 \hline
 \end{tabular}
 \end{center}
 \vskip0.3cm
 It then follows that
 $X_{2n}(t)|_{U_2}-\frac 14 X_n(t^{(2)})$ has poles only at
 $\infty$ with $q^{-n}/(4n)$ as its pole part. Thus
 \begin{equation}
 X_{2n}(t)|_{U_2}-\frac 14 X_n(t^{(2)})=\frac 14 X_n(t).
 \label{IIa}
 \end{equation}
 And $X_{2n}(t_0)|_{U_2}-\frac 14 X_n(t_0^{(2)})$ has poles only at
 $\infty$ and $3/10$ with $q^{-n}/(4n)$ as its pole parts.
 Therefore,
 \begin{equation}
 X_{2n}(t_0)|_{U_2}-\frac 14 X_n(t_0^{(2)})=\frac 14 X_n(t_0).
 \label{IIb}
 \end{equation}
 Now subtracting (\ref{IIb}) from two times the equality (\ref{IIa})
 leads us to the conclusion.
 \end{pf}
 \vskip0.1cm
 \noindent Let $\chi_0$ be the trivial character mod $2$, i.e.
 $\chi_0 (n)=
  \begin{cases} 0, & \text{ if $n$ is even } \\
                1, & \text{ if $n$ is odd. }
  \end{cases} $
 \par\noindent For a modular function $f=\sum_{n\in {\Bbb Z}} a_n q^n$, we
 define its {\it twist} by $f_{\chi_0}=\sum_{n\in {\Bbb Z}} a_n \chi_0 (n)
 q^n$. Then it is easy to get
 $f_{\chi_0}=f-f|_{U_2}(2z)=\frac 12 (f-f|_{\sm 2&1\\0&2 \esm}).$
 \begin{lem}
 For an odd positive integer $r$,
 \begin{equation}
 [(2X_{r}(t)- X_{r}(t_0))|_{U_{2^k}}]_{\chi_0}
 =(-1)^k [2X_{2^k r}(t)- X_{2^k r}(t_0)]_{\chi_0}.
 \label{JJa}
 \end{equation}
 \label{JJ}
 \end{lem}
 \begin{pf}
 Since $t_0$ is replicable, it follows that for odd $m$, $h_{2^k m, r}=h_{m, 2^k
 r}$, whence $[X_r(t_0)|_{U_{2^k}}]_{\chi_0}=X_{2^k
 r}(t_0)_{\chi_0}$. Thus
 (\ref{JJa}) is equivalent to
 $$ [2X_r(t)|_{U_{2^k}}]_{\chi_0}-X_{2^k r}(t_0)_{\chi_0}
    =(-1)^k \cdot 2\cdot X_{2^k r}(t)_{\chi_0}+
    (-1)^{k+1} X_{2^k r}(t_0)_{\chi_0}.$$
 Now to justify (\ref{JJa}) it suffices to prove the following.
 \begin{equation}
 [X_r(t)|_{U_{2^k}}]_{\chi_0}=(-1)^k \cdot X_{2^k r}(t)_{\chi_0}
   + \frac{(-1)^{k+1}+1}{2}\cdot X_{2^k r}(t_0)_{\chi_0}.
 \label{JJb}
 \end{equation}
 Since $f_{\chi_0}=f-f|_{U_2}(2z)$, both LHS and RHS of
 (\ref{JJb}) belong to $K(X_1(20))$.
 \begin{align*}
 \text{ LHS of (\ref{JJb}) }
 & = X_r(t)|_{U_{2^k}} - X_r(t)|_{U_{2^{k+1}}}(2z)=
     \frac 12
     \left[X_r(t)|_{U_{2^k}}-X_r(t)|_{U_{2^k}}(z+1/2)\right] \\
 & = \frac{1}{2^{k+1}}\sum_{i=0}^{1}\sum_{j=0}^{2^k -1}
     X_r(t)|_{\sm 1&j \\ 0&2^k \esm \sm 2&i\\0&2 \esm},
 \end{align*}
 which can have poles only at $\sm 2&i \\ 0&2 \esm^{-1} \sm 1&j\\0&2^k
 \esm^{-1}\Gamma_1(10)\infty = \sm 2^{k+1}&-2j-i \\ 0&2 \esm
 \Gamma_1(10)\infty$.
 Let $\sm a&b\\c&d \esm\in \Gamma_1(10)$. Then
 $\sm 2^{k+1}&-2j-i \\ 0&2 \esm
 \sm a&b\\c&d \esm\infty=(2^{k+1} a -2jc-ic)/(2c).$
 Observe that g.c.d. of $2^{k+1} a -2jc-ic$ and $2c$ divides
 $2^{k+2}$. Write it as $2^l$.
 Then $s=\frac{(2^{k+1} a -2jc-ic)/2^l}{2c/2^l}$ is of the form
 $n/(5m)$ for some integer $m$ and $n$ with $(5,n)=1$.
 Similarly, RHS of (\ref{JJb}) can have poles only at such cusps.
 Upto $\Gamma_1(20)$-equivalence, cusp of the form $n/(5m)$ is
 equivalent to one of $1/5,2/5,3/5,4/5, 1/10, 3/10, 1/20, 3/20,
 7/20, 9/20$.
 \par\noindent
 \underline{$s=1/20, 3/20, 7/20, 9/20$}
 \par\noindent
 We claim that for $f\in K(X_1(10))$ and
 $\gamma_0=\sm a&b\\c&d \esm\in \Gamma_0(20)$,
 \begin{equation}
 (f_{\chi_0})|_{\gamma_0}=
 (f|_{\gamma_0})_{\chi_0}.
 \label{JJc}
 \end{equation}
 Indeed, the LHS of the above is $(f-f|_{U_2 \sm
 2&0\\0&1\esm})|_{\gamma_0}=f|_{\gamma_0}-f|_{U_2 \sm
 2&0\\0&1\esm\gamma_0}$. Meanwhile, the RHS of (\ref{JJc}) is
 $f|_{\gamma_0}-(f_{\gamma_0})|_{U_2 \sm
 2&0\\0&1\esm}=f|_{\gamma_0}-(f|_{U_2})|_{\gamma_0 \sm
 2&0\\0&1\esm}$ by Lemma \ref{Hecke}.
 In order to claim (\ref{JJc}) it is
 enough to show that $\gamma_0 \sm 2&0\\0&1 \esm {\gamma_0}^{-1}
 \sm 2&0\\0&1 \esm^{-1} \in \Gamma_1(10)$. In fact,
 $\sm a&b\\c&d \esm \sm 2&0\\0&1 \esm {\sm a&b\\c&d \esm}^{-1}
 \sm 2&0\\0&1 \esm^{-1}=\sm ad-bc/2& -ab\\ cd/2& -2bc+da
 \esm\in\Gamma_1(10)$.
 Now we take $f=X_r(t)|_{U_{2^k}}$ (or $X_{2^k r}(t)$ or $X_{2^k
 r}(t_0)$) and $\gamma_0\in\Gamma_0(20)$ which sends $\infty$ to $s$.
 Then
 \begin{align*}
 [\text{ LHS of (\ref{JJb}) }]|_{\gamma_0}
 & =[X_r(t)|_{U_{2^k}}]_{\chi_0}|_{\gamma_0}
   =[[X_r(t)|_{U_{2^k}}]_{\gamma_0}]_{\chi_0} \\
 & =[X_r(t)|_{\gamma_0 U_{2^k}}]_{\chi_0} \text{ by
 Lemma \ref{Hecke} }
 \end{align*}
 which belongs to $O(1)$. And
 \par\noindent
 $[\text{ RHS of (\ref{JJb}) }]|_{\gamma_0}
 =(-1)^k X_{2^k r}(t)_{\chi_0}|_{\gamma_0}
  +\frac{(-1)^{k+1}+1}{2}\cdot X_{2^k
  r}(t_0)_{\chi_0}|_{\gamma_0}$. Here we observe that
 $$X_{2^k r}(t)_{\chi_0}|_{\gamma_0}=[X_{2^k r}(t)|_{\gamma_0}]_{\chi_0}
  =\begin{cases} X_{2^k r}(t)_{\chi_0}, & \text{ if } \gamma_0\in\pm
  \Gamma_1(10) \\ O(1), & \text{ otherwise } \end{cases}$$
  which belongs to $O(1)$.
 In a similar way, we get $X_{2^k r}(t_0)_{\chi_0}|_{\gamma_0}\in O(1)$.
 Thus $[\text{ RHS of (\ref{JJb}) }]|_{\gamma_0}\in O(1)$.
 Now LHS of (\ref{JJb}) $-$ RHS of (\ref{JJb}) is holomorphic at
 the cusps
 $s=1/20,3/20, 7/20,9/20$.
 \par\noindent
 To deal with other cusps we need a sublemma.
 \par\noindent
 {\bf Sublemma.}
 \begin{align*}
 [X_r(t)|_{U_{2^k}}]_{\sm 1&0\\5&1 \esm}
 &= \begin{cases} \frac{1}{2r} q_2^{-r}, & \text{ if } k=1 \\
    \frac{(-1)^{k+1}+1}{2}\cdot \left(-\frac{1}{2^k r}\right)
    \cdot q^{-2^{k-2} r} + O(1), & \text{ if } k\ge 2
    \end{cases} \\
 \intertext{ and }
 [X_r(t_0)|_{U_{2^k}}]_{\sm 1&0\\5&1 \esm}
 &= \begin{cases} \frac{1}{2r} q_2^{-r}, & \text{ if } k=1 \\
    -\frac{1}{2^k r}\cdot q^{-2^{k-2} r} + O(1), & \text{ if } k\ge
    2.
    \end{cases}
 \end{align*}
 \par\noindent
 (proof of the sublemma)
 \begin{align*}
 [X_r(t)|_{U_{2^k}}]_{\sm 1&0\\5&1 \esm}
 & = [X_r(t)|_{U_{2^{k-1}}}]_{U_2 \sm 1&0\\5&1 \esm}
   =\frac 12 [X_r(t)|_{U_{2^{k-1}}}]_{\sm 1&0\\0&2 \esm \sm 1&0\\5&1 \esm}
   +\frac 12 [X_r(t)|_{U_{2^{k-1}}}]_{\sm 1&1\\0&2 \esm \sm 1&0\\5&1
   \esm}\\
 & =\frac 12 [X_r(t)|_{U_{2^{k-1}}}]_{\sm 1&0\\10&1 \esm \sm 1&0\\0&2 \esm}
   +\frac 12 [X_r(t)|_{U_{2^{k-1}}}]_{\sm 13&-2\\20&-3 \esm \sm 1&0\\5&1
   \esm \sm 2&1\\0&1 \esm} \\
 & =\frac 12 X_r(t)|_{U_{2^{k-1}}}(z/2)+\frac
 12 [X_r(t_0)-X_r(t)]|_{U_{2^{k-1}}\sm 1&0\\5&1
   \esm \sm 2&1\\0&1 \esm} + \text{ const. } \\
 & \text{ \q \q \q by Lemma \ref{Hecke} and Lemma \ref{Action} }.
 \end{align*}
 Thus
 \begin{equation}
 [X_r(t)|_{U_{2^k}}]_{\sm 1&0\\5&1 \esm}
 = \begin{cases} \frac{1}{2r} q_2^{-r}, & \text{ if } k=1 \\
    \frac{1}{2}[X_r(t_0)-X_r(t)]|_{U_{2^{k-1}}\sm 1&0\\5&1
   \esm \sm 2&1\\0&1 \esm} + O(1), & \text{ if } k\ge
    2.
   \end{cases}
 \label{JJd}
 \end{equation}
 Likewise, we have
 \begin{equation}
 [X_r(t_0)|_{U_{2^k}}]_{\sm 1&0\\5&1 \esm}
 = \begin{cases} \frac{1}{2r} q_2^{-r}, & \text{ if } k=1 \\
    \frac{1}{2}X_r(t_0)|_{U_{2^{k-1}}\sm 1&0\\5&1
   \esm \sm 2&1\\0&1 \esm} + O(1), & \text{ if } k\ge
    2.
   \end{cases}
 \label{JJe}
 \end{equation}
 Now we prove the sublemma by induction on $k$.
 \par\noindent
 If $k=1$, it was done.
 \par\noindent
 If $k=2$,
 \begin{align*}
 X_r(t)|_{U_4\sm 1&0\\5&1 \esm}
 & = \frac{1}{2}[X_r(t_0)-X_r(t)]|_{U_2 \sm 1&0\\5&1
     \esm \sm 2&1\\0&1 \esm} + O(1) \\
 & = O(1) \text{ by (\ref{JJd}) and (\ref{JJe}) } \\
 \intertext{ and }
 X_r(t_0)|_{U_4\sm 1&0\\5&1 \esm}
 & = \frac{1}{2} X_r(t_0)_{U_2 \sm 1&0\\5&1
     \esm \sm 2&1\\0&1 \esm} + O(1) = \frac{1}{4r}\cdot e^{-\pi ir(2z+1)}+O(1) \\
 & = -\frac{1}{4r} q^{-r} + O(1).
 \end{align*}
 Now let $k\ge 3$.
 \begin{align*}
 [X_r(t)|_{U_{2^k}}]_{\sm 1&0\\5&1 \esm}
 & = \frac 12 X_r(t_0)|_{U_{2^{k-1}}\sm 1&0\\5&1 \esm}(2z+1)
    -\frac 12 X_r(t)|_{U_{2^{k-1}}\sm 1&0\\5&1 \esm}(2z+1) + O(1)
    \text{ by (\ref{JJd}) } \\
 & = \frac 12\left(-\frac{1}{2^{k-1}r} e^{-2\pi i\cdot 2^{k-3}r(2z+1)}\right)
    -\frac 12\cdot\frac{(-1)^k+1}{2}\cdot\left(-\frac{1}{2^{k-1}r}\right)
    e^{-2\pi i\cdot 2^{k-3}r(2z+1)}+O(1) \\
 & \text{ \q \q \q by induction hypothesis for $k-1$ } \\
 & =-\frac{1}{2^{k}r} q^{-2^{k-2} r}
    +\frac{(-1)^k+1}{2}\cdot \frac{1}{2^k r} q^{-2^{k-2} r}+O(1)
    \\
 & =-\frac{1}{2^{k}r} q^{-2^{k-2} r}\left(
   1-\frac{(-1)^k+1}{2}\right)+O(1) \\
 & =-\frac{1}{2^{k}r} q^{-2^{k-2} r}\left(
   \frac{1+(-1)^{k+1}}{2}\right)+O(1).
 \end{align*}
 In like manner, the assertion for $[X_r(t_0)|_{U_{2^k}}]_{\sm 1&0\\5&1
 \esm}$ can be proved. This completes the proof of the sublemma.
 \par\noindent
 \underline{$s=1/5$}
 \par\noindent
 \begin{align*}
 [\text{ LHS of (\ref{JJb}) }]|_{\sm 1&0\\5&1 \esm}
 & =[X_r(t)|_{U_{2^k}}]_{\chi_0}|_{\sm 1&0\\5&1 \esm}
   =X_r(t)|_{U_{2^k}\sm 1&0\\5&1 \esm}
    -X_r(t)|_{U_{2^{k+1}}\sm 2&0\\0&1 \esm\sm 1&0\\5&1 \esm} \\
 & = X_r(t)|_{U_{2^k}\sm 1&0\\5&1 \esm}
    -X_r(t)|_{U_{2^{k+1}}\sm -3&1\\-10&3 \esm\sm 1&0\\5&1 \esm\sm 1&-1\\0&2
    \esm}\\
 & = X_r(t)|_{U_{2^k}\sm 1&0\\5&1 \esm}
    -[X_r(t_0)-X_r(t)]|_{U_{2^{k+1}}\sm 1&0\\5&1 \esm\sm 1&-1\\0&2
    \esm}+\text{ const. } \\
 & \text{ \q \q \q by Lemma \ref{Hecke} and Lemma \ref{Action} } \\
 & = \begin{cases}
     \frac{1}{2r} q_2^{-r}+\frac{1}{4r}e^{-2\pi i
     r(\frac{z-1}{2})}+O(1), & \text{ if } k=1 \\
     O(1), & \text{ if $k$ even } \\
     -\frac{1}{2^k r}q^{-2^{k-2}r}
     +\frac{1}{2^{k+1}r}q^{-2^{k-2}r}+O(1), & \text{ if $k\ge3$
     odd }
     \end{cases}
     \text{ by the sublemma } \\
 & = \begin{cases}
     \frac{1}{4r} q_2^{-r}, & \text{ if } k=1 \\
     \frac{(-1)^{k+1}+1}{2}\cdot \left(-\frac{1}{2^{k+1} r}
     \cdot q^{-2^{k-2} r}\right) + O(1), & \text{ if } k\ge 2.
     \end{cases}
 \end{align*}
 Observe that $f_{\chi_0}|_{\sm 1&0\\5&1 \esm}=\frac 12 [f-f|_{\sm
 2&1\\0&2\esm}]|_{\sm 1&0\\5&1 \esm}=\frac 12 f|_{\sm 1&0\\5&1 \esm}
 -\frac 12 f|_{\sm 7&2\\10&3 \esm\sm 1&-1\\0&4 \esm}$.
 If we take $X_{2^k r}(t)$ (and $X_{2^k r}(t_0)$) instead of $f$, it then follows that
 $X_{2^k r}(t)_{\chi_0}|_{\sm 1&0\\5&1 \esm}=
  \frac 12 X_{2^k r}(t)|_{\sm 1&0\\5&1 \esm}
 -\frac 12 X_{2^k r}(t)|_{\sm 7&2\\10&3 \esm\sm 1&-1\\0&4 \esm}\in
 O(1)$ and
 \begin{align*}
 X_{2^k r}(t_0)_{\chi_0}|_{\sm 1&0\\5&1 \esm}
 & = \frac 12 X_{2^k r}(t_0)|_{\sm 1&0\\5&1 \esm}
 -\frac 12 X_{2^k r}(t_0)|_{\sm 7&2\\10&3 \esm\sm 1&-1\\0&4 \esm}=
 -\frac 12 \cdot \frac{1}{2^k r}\cdot e^{-2\pi i 2^k
 r(\frac{z-1}{4})}+O(1) \\
 & =\begin{cases}
    \frac{1}{4r}q_2^{-r}+O(1), & \text{ if } k=1 \\
    -\frac{1}{2^{k+1} r}\cdot q^{-2^{k-2} r} + O(1), & \text{ if } k\ge 2.
     \end{cases}
 \end{align*}
 Therefore
 $$[\text{ RHS of (\ref{JJb}) }]|_{\sm 1&0\\5&1 \esm}
  =\begin{cases}
   \frac{1}{4r}q_2^{-r}+O(1), & \text{ if } k=1 \\
   \frac{(-1)^{k+1}+1}{2}\cdot \left(-\frac{1}{2^{k+1} r}
     \cdot q^{-2^{k-2} r}\right) + O(1), & \text{ if } k\ge 2.
     \end{cases} $$
 Now LHS $-$ RHS is holomorphic at $s=1/5$.
 \par\noindent
 \underline{$s=1/10$}
 For all $k\ge 1$,
 \begin{align*}
 [\text{ LHS of (\ref{JJb}) }]|_{\sm 1&0\\10&1 \esm}
 & =[X_r(t)|_{U_{2^k}}]_{\chi_0}|_{\sm 1&0\\10&1 \esm}
   =X_r(t)|_{U_{2^k}\sm 1&0\\10&1 \esm}
   -X_r(t)|_{U_{2^{k+1}}\sm 2&0\\0&1 \esm\sm 1&0\\10&1 \esm} \\
 & = X_r(t)|_{U_{2^k}}-X_r(t)|_{U_{2^{k+1}}\sm 1&0\\5&1 \esm\sm 2&0\\0&1 \esm} \\
 & = \frac{(-1)^{k}+1}{2}\cdot \frac{1}{2^{k+1} r}
     \cdot q^{-2^{k} r} + O(1) \text{ by the sublemma. }
 \end{align*}
 Observe that $f_{\chi_0}|_{\sm 1&0\\10&1 \esm}=\frac 12 [f-f|_{\sm
 2&1\\0&2\esm}]|_{\sm 1&0\\10&1 \esm}=\frac 12 f|_{\sm 1&0\\10&1 \esm}
 -\frac 12 f|_{\sm 3&1\\5&2 \esm\sm 4&0\\0&1 \esm}$.
 If we take $f$ to be $X_{2^k r}(t)$ (or $X_{2^k r}(t)$), then it follows that
 \begin{align*}
 X_{2^k r}(t)_{\chi_0}|_{\sm 1&0\\10&1 \esm}
 & =\frac 12 X_{2^k r}(t)|_{\sm 1&0\\10&1 \esm}
    -\frac 12 X_{2^k r}(t)|_{\sm 3&1\\5&2 \esm\sm 4&0\\0&1 \esm}
    \\
 & = \frac{1}{2^{k+1}r} q^{-2^k r} + O(1) \\
 \intertext{ and }
 & X_{2^k r}(t_0)_{\chi_0}|_{\sm 1&0\\10&1 \esm}
  = \frac{1}{2^{k+1}r} q^{-2^k r} + O(1).
 \end{align*}
 Thus
 \begin{align*}
 [\text{ RHS of (\ref{JJb}) }]|_{\sm 1&0\\10&1 \esm}
 & = \left( (-1)^k + \frac{(-1)^{k+1}+1}{2}\right)\cdot \frac{1}{2^{k+1} r}
    \cdot q^{-2^{k} r} + O(1) \\
 & = \frac{(-1)^{k}+1}{2}\frac{1}{2^{k+1} r} \cdot q^{-2^{k} r} +
 O(1).
 \end{align*}
 Now LHS $-$ RHS is holomorphic at $s=1/10$.
 \par\noindent
 \underline{$s=2/5,3/5,4/5,3/10$}
 \par\noindent
 Observe that $2/5\sim 3/5 \sim 4/5 \sim 1/5$ and $3/10\sim 1/10$
 under $\Gamma_0(20)$. Indeed, for $c|N$, $a/c\sim a'/c$ in $X_0(20)$ if and
 only if $a\equiv a' \mod (c,N/c)$.
 Write $s=\gamma_0 s_0$ with $\gamma_0\in \Gamma_0(20)$ and
 $s_0=1/5$ or $1/10$.
 In (\ref{JJc}) we take $f$ to be
 $[2X_r(t)-X_r(t_0)]|_{U_{2^k}}-(-1)^k\cdot [2X_{2^k r}(t)-X_{2^k
 r}(t_0)]$. By (\ref{JJc}), Lemma \ref{Hecke} and Lemma \ref{Action}, we obtain that
 $$(f_{\chi_0})|_{\gamma_0}=(f|_{\gamma_0})_{\chi_0}=
  \begin{cases}
  f_{\chi_0}, & \text{ if } \gamma_0\in\pm \Gamma_1(10) \\
  (-f)_{\chi_0} + \text{ constant, } & \text{ if } \gamma_0\not\in\pm \Gamma_1(10).
  \end{cases} $$
 Note that $\gamma_0\sm 1&0\\5&1 \esm$ or $\gamma_0\sm 1&0\\10&1
 \esm$ sends $\infty$ to $s$ according as $s_0=1/5$ or $s_0=1/10$.
 Moreover $f_{\chi_0}$ is equal to LHS of (\ref{JJa}) $-$
 RHS of (\ref{JJa}).
 Finally
 $ (f_{\chi_0})|_{\gamma_0 \sm 1&0\\5&1 \esm}
 = \pm (f_{\chi_0})|_{\sm 1&0\\5&1 \esm} \in O(1). $ And
 $ (f_{\chi_0})|_{\gamma_0\sm 1&0\\10&1 \esm}
 = \pm (f_{\chi_0})|_{\sm 1&0\\10&1 \esm} \in O(1). $
 \end{pf}
 \vskip0.2cm
 For positive integers $a$ and $b$ we define $\psi (a,b)$ as follows.
 First we assume that $(10,(a,b))=1$.
 \par\noindent
 Case I: $a$ odd and $b$ odd.
 $$\psi (a,b) =\begin{cases} 1, & \text{ if } a\equiv \pm 1
 \mod 10 \text{ or } b\equiv \pm 1 \mod 10 \\
 -1, & \text{ otherwise.} \end{cases}$$
 \par\noindent
 Case II: $a$ even and $b$ odd. Write $a=2^k r$ with
 $r$ odd.
 $$\psi (a,b) =\begin{cases} \psi_{10} (b), & \text{ if } (10,b)=1 \\
 (-1)^k\cdot \psi_{10} (r), & \text{ otherwise.} \end{cases}$$
 \par\noindent
 Case III: $a$ odd and $b$ even.
 \par Put $\psi (a,b)= \psi (b,a)$.
 \par\noindent
 In general, let $2^\alpha 5^\beta || (a,b)$ and define
 $\psi (a,b) = \psi (\frac{a}{2^\alpha 5^\beta},\frac{b}{2^\alpha
 5^\beta})$.
 \par\noindent
 We write $X_n(t)=\sum_{m\in {\Bbb Z}} H_{m,n} q^m$,
 $X_n(t_0)=\sum_{m\in {\Bbb Z}} h_{m,n} q^m$,
 $X_n(t^{(2)})=\sum_{m\in {\Bbb Z}} H_{m,n}^{(2)} q^m$, and
 $X_n(t_0^{(2)})=\sum_{m\in {\Bbb Z}} h_{m,n}^{(2)} q^m$.
 \begin{theo}
 Let $N=10$. If $ab=cd$ and $(a,b)=(c,d)$, then
 $$ \psi (a,b) \times (2H_{a,b}-h_{a,b}) =
 \psi (c,d) \times (2H_{c,d}-h_{c,d}).$$
 \label{KK}
 \end{theo}
 \begin{pf}
 Let $2^\alpha 5^\beta || (a,b)=(c,d)$.
 Write $a'=\frac{a}{2^\alpha 5^\beta}$,
 $b'=\frac{b}{2^\alpha 5^\beta}$,
 $c'=\frac{c}{2^\alpha 5^\beta}$, and
 $d'=\frac{d}{2^\alpha 5^\beta}$.
 By Lemma \ref{AA} and Lemma \ref{II},
 \begin{align}
   2H_{a,b}-h_{a,b} & =\frac{1}{5^\beta}
    (2H_{\frac{a}{5^\beta},\frac{b}{5^\beta}}
     -h_{\frac{a}{5^\beta},\frac{b}{5^\beta}})
     \label{Ka} \\
     & = \frac{1}{4^\alpha 5^\beta}
         (2H_{a',b'}^{(2)}-h_{a',b'}^{(2)}
          +2H_{a',b'}-h_{a',b'}). \notag
 \end{align}
 Similarly,
 \begin{equation}
 2H_{c,d}-h_{c,d}=\frac{1}{4^\alpha 5^\beta}
    (2H_{c',d'}^{(2)}-h_{c',d'}^{(2)}+2H_{c',d'}-h_{c',d'}).
 \label{Kb}
 \end{equation}
 \par\noindent
 Case I: $a'$ odd and $b'$ odd
 \par
 In this case $(a',10)=1$ or $(b',10)=1$. And $(c',10)=1$ or
 $(d',10)=1$. By $(*)$
 it holds that $\psi (a',b') \times (\ref{Ka})
 = \psi (c',d') \times (\ref{Kb})$. This implies the assertion.
 \par\noindent
 Case II: $a'$ even or $b'$ even
 \par
 Since $\psi (a',b')=\psi (b',a')$, $H_{a',b'}=H_{b',a'}$ and
 $h_{a',b'}=h_{b',a'}$, we may assume that $a'$ even and $b'$ odd
 (resp. $c'$ even and $d'$ odd).
 Write $a'=2^k\cdot r$ with $r$ odd (resp.
 $c'=2^k\cdot s$ with $s$ odd) and $k\ge 1$.
 Then $(r,10)=1$ or $(b',10)=1$. And $(s,10)=1$ or $(d',10)=1$.
 \par\noindent
 If $(b',10)=(d',10)=1$, then by $(*)$
 \begin{align*}
 & \psi_5 (b') \times (2H_{a',b'}^{(2)}-h_{a',b'}^{(2)}) =
 \psi_5 (d') \times (2H_{c',d'}^{(2)}-h_{c',d'}^{(2)})  \\
 \intertext{ and }
 & \psi_{10} (b') \times (2H_{a',b'}-h_{a',b'}) =
 \psi_{10} (d') \times (2H_{c',d'}-h_{c',d'}).
 \end{align*}
 Notice that $\psi_5 (n)= \psi_{10} (n)$ for $n$ prime to $10$.
 We then have $\psi (a,b) \times (\ref{Ka})
 = \psi (c,d) \times (\ref{Kb})$.
 \par\noindent
 If $(b',10)\neq 1$ and $(d',10)=1$, then $(r,10)=1$.
 Note that $(a',5)=1$ and $(d',5)=1$. By
 $(*)$
 $$
  \psi_5 (a') \times (2H_{a',b'}^{(2)}-h_{a',b'}^{(2)}) =
 \psi_5 (d') \times (2H_{c',d'}^{(2)}-h_{c',d'}^{(2)}).$$
 Since $\psi_5 (2^k r)=\psi_5(2)^k \psi_5 (r)=(-1)^k \psi_{10} (r)$,
 the above can be rewritten as
 $$ (-1)^k \psi_{10} (r)\times (2H_{a',b'}^{(2)}-h_{a',b'}^{(2)}) =
 \psi_{10} (d') \times (2H_{c',d'}^{(2)}-h_{c',d'}^{(2)}).$$
 By Lemma \ref{JJ},
 \begin{align}
 (-1)^k (2H_{a',b'}-h_{a',b'})
 & = (-1)^k (2H_{2^k r,b'}-h_{2^k r,b'})=(-1)^k (2H_{b',2^k r}-h_{b',2^k r})
  \label{Kc} \\
 & = 2H_{2^k b',r}-h_{2^k b',r}.
  \notag
 \end{align}
 Since $(r,10)=(d',10)=1$, we obtain by $(*)$
 that $\psi_{10} (r) \times (\ref{Kc})
 = \psi_{10} (d') \times (2H_{c',d'}-h_{c',d'})$. This implies
 $(-1)^k \psi_{10} (r) \times (2H_{a',b'}-h_{a',b'})
  = \psi_{10} (d') \times (2H_{c',d'}-h_{c',d'})$. Thus
 $\psi (a,b) \times (\ref{Ka})
 = \psi (c,d) \times (\ref{Kb})$.
 The remaining cases (i.e. $(b',10)=1$ and $(d',10)\neq 1$
 or $(b',10)\neq 1$ and $(d',10)\neq 1$) can be treated in a similar way.
 \end{pf}
 
\end{document}